# Quantum optimisation in cities: Limitations and prospects of urban transport systems

Junxiang Xu[1,*], Chence Niu[2], Divya Jayakumar Nair[1], Vinayak Dixit[1]

1. Research Centre for Integrated Transport Innovation (rCITI), School of Civil and Environmental Engineering, The University of New South Wales, Kensington, UNSW Sydney, NSW, 2052, Australia
2. Guangdong Basic Research Center of Excellence for Ecological Security and Green Development, Key Laboratory for City Cluster Environmental Safety and Green Development of the Ministry of Education, School of Ecology, Environment and Resources, Guangdong University of Technology, Guangzhou, 510006, China

**corresponding author\*:** junxiang.xu@unsw.edu.au

**Abstract:** Recently, quantum computing has attracted extensive attention in urban studies as a potential tool for addressing complex decision problems in urban transport systems. However, its role in transport planning remains poorly defined. This paper reviews the existing research on quantum computing in urban transport planning and points out the substantial limitations of current quantum optimisation methods in terms of scalability, robustness, constraint representation, and engineering feasibility. So, the stable and reproducible advantages of quantum optimisation in real urban systems remain to be demonstrated. By comparing quantum methods with established classical optimisation methods, it is found that decomposition methods, metaheuristics, and reinforcement learning already provide transparent, scalable, and policy-interpretable solutions for medium and large-sized urban transport networks. In contrast, the contribution of quantum methods largely lies in the exploratory analysis of limited, discrete combinatorial subproblems rather than full system-level optimisation. It is argued in this paper for a shift from technology-driven application narrative towards problem-driven method selection. From an urban transport planning perspective, we have identified the specific problem types where the exploratory use of quantum computing may be relevant, including critical link and node vulnerability identification, combinatorial screening of congestion and failure scenarios, disaster-related extreme condition analysis, constrained path option selection, and small-scale facility location and investment option assessment. It is concluded that hybrid frameworks represent a more realistic pathway for integrating quantum computing into urban transport research, in which classical methods ensure system-level consistency and policy interpretability while quantum methods support local combinatorial exploration. Until stable engineering advantages are demonstrated, public agencies and researchers should prioritise method validation, scenario suitability, and cross-disciplinary collaboration.

**Keywords:** quantum computing; urban transport planning; combinatorial optimisation; hybrid decision framework; urban policy implication

## 1. Introduction

*1.1 Urban transport systems as complex decision environments*

Urban transport systems form complex public decision systems. Their planning and operation involve strategic,

tactical, and operational decision levels (Cascetta et al., 2015). These levels differ clearly in time, space, and management scope, yet they remain coupled throughout the transport network structure and operating states (Lin and Ban, 2013). As a result, urban transport problems cannot be reduced to single-level technical optimisation tasks. From a structural perspective, urban transport systems are typically featured by large network scale, high-dimensional decision variables, coexistence of continuous and discrete decisions, and overlapping operational, regulatory, and institutional constraints (Pavez et al., 2025; Zou et al., 2025). These features determine the complexity in modelling and solving the transport system, and place high demands on the applicability and practical feasibility of computational methods (Chen et al., 2011; Barbati et al., 2012; Bastarianto et al., 2023). As illustrated in **Figure 1**, urban transport decision-making can be categorised into strategic, tactical, and operational levels. Each level corresponds to distinct characteristics of transport problems, and together they define the decision environment for selecting the appropriate computational and methodological toolboxes.

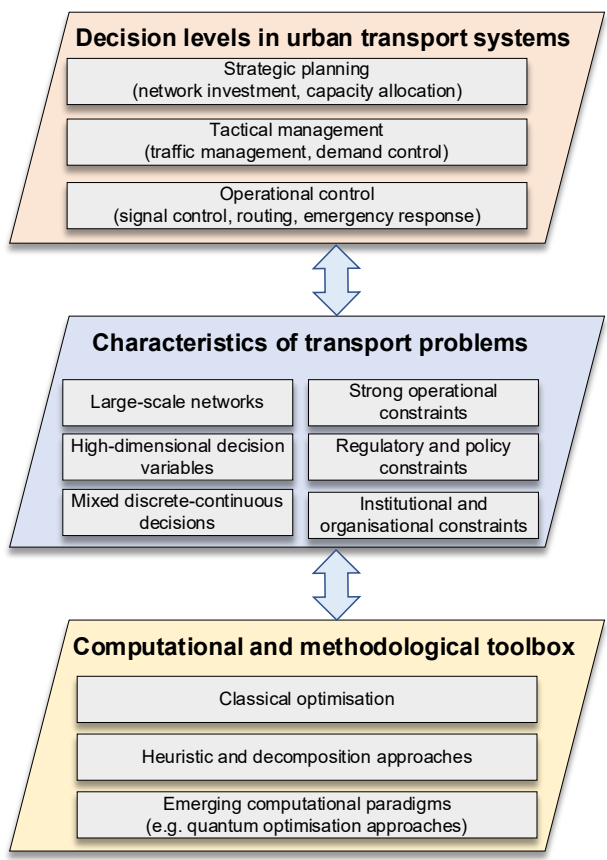

**Figure 1.** Urban transport systems as a multi-level decision environment.

*1.2 Existing research on quantum computing in urban transport planning*

Although quantum computing has received growing attention in combinatorial optimisation in recent years, its application in urban transport systems remains at an early stage, with limited literature and relatively concentrated research directions. Existing studies mainly attempt to apply typical problems in urban transport to quantum computing frameworks. For example, Cooper (2021) innovatively explored the potential future applications of quantum computing in transport simulation and planning. Dixit and Jian (2022) used quantum technology to achieve the fast estimation of vehicle driving cycle frequency. Their work showed great importance for energy saving, emission reduction, and safety improvement. They later further extended quantum computing to transport network design and laid an initial foundation for introducing quantum technology into urban transport planning (Dixit and Niu, 2023). Qu et al. (2022) creatively integrated quantum computing with machine learning and developed a spatiotemporal quantum convolutional neural network for traffic congestion prediction. Schetakis et al. (2025) further demonstrated the reliability and practical value of quantum neural networks in traffic congestion prediction. This approach provides technical support for the further improvement of urban short-term traffic flow control. Dixit et al. (2024) proposed a quantum optimisation approach to search for optimal paths using real-time-dependent traffic data. In addition, quantum-inspired heuristics developed by Tian et al. (2019) were validated for path planning in large-scale datasets. Liu et al. (2025) pointed out the concerns about the extensive application of quantum technology in current transport and logistics fields, which aligns with the motivation of this study. However, this study focuses on urban transport planning. It supplements the details on quantum computing models and application scenarios and identifies feasible directions for future implementation. Recent studies applied quantum annealing (QA) to optimise urban multi-commodity flow problems and compared computational advantages of QA algorithms over heuristic methods (Niu et al., 2025b). Du et al. (2026) developed a framework to assess the potential scalability of quantum algorithms in air transport network design and considered the influences of multiple factors, such as qubit numbers and error rates. However, they argued that real-world application of quantum technology remained distant, mainly based on reflections on the limitations of quantum hardware and algorithms. In summary, it can be observed that the existing research linking quantum technologies with transport planning mainly focuses on transport network design, traffic signal control, route planning, traffic congestion prediction, and intelligent transport system design. **Figure 2** presents the research timeline of quantum computing applications and development.

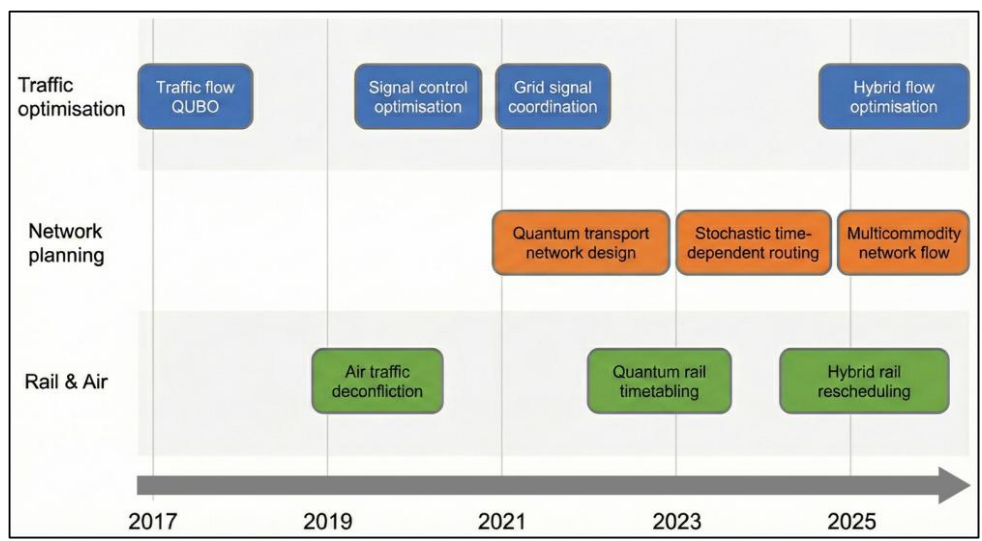

**Figure 2.** Research timeline of quantum computing in transport studies 2017 to 2025 (The supporting references in this figure are provided in **Appendix A**).

Overall, the existing studies have largely focused on several specific topics, e.g. problem mapping methods, variable encoding strategies, and computational performance of quantum or quantum-inspired solution frameworks under benchmark networks or case scenarios. Their research objectives mainly concentrate on method feasibility verification or performance comparison with classical algorithms. In contrast, there is rare discussion on the application scope of quantum optimisation methods across different decision levels in urban transport systems and their capacity to extend scientific questions. Although some of the above studies have begun to consider the limitations and generalisability, they did not provide constructive insights on how to practically solve the challenges in quantum computing-based urban transport planning.

*1.3 Critical reflections on existing quantum optimisation studies*
Existing studies have noticed that the practical deployment of quantum optimisation still faces substantive obstacles. First, a consistent and robust demonstration of *quantum advantage* on large-scale, real-world instances has yet to be achieved. Many reported successes rely on carefully-selected synthetic or toy benchmarks, where the design of test cases and the criteria for comparison remain open to bias and contestation (Abbas et al., 2024). At the same time, the challenges from scalability, noise and robustness across both algorithms and hardware have limited the prospects of systematic breakthroughs. For example, recent assessments of variational quantum algorithms indicate that their scalability and superiority over strong classical baselines remain far from being conclusive (Müller et al., 2025). Moreover, field-specific reviews in transport suggest that the current experiments with Quadratic Unconstrained Binary Optimisation (QUBO) reformulations on quantum annealers or hybrid solvers largely stay at the level of problem compatibility and benchmark exploration, and that consistent outperformance on realistic network scales remains elusive (Mohammed et al., 2025).

Currently, several exploratory studies in the transport field have reformulated problems such as the vehicle routing problem (VRP), facility location-allocation, and evacuation route planning into QUBO form, enabling their execution on quantum annealers or hybrid quantum-classical solvers (Leonidas et al., 2023; Cattelan and Yarkoni, 2024). These efforts claim that quantum computing can provide new ways to address complex scheduling and routing tasks. However, further investigation reveals that these problems can already be solved by a wide range of mature classical methods. Mixed Integer Linear Programming (MILP) has been a common and effective method to deal with transport cases (Luathep et al., 2011; Demirel et al., 2016; Yuan et al., 2019; Baller et al., 2022). Decomposition strategies, such as Benders decomposition and column generation, can address large-scale formulations (Zeighami and Soumis, 2019; Lan et al., 2021; Li et al., 2022b). Heuristics and metaheuristics (e.g. genetic algorithms, simulated annealing, tabu search) have been extensively employed in industrial practice (Xu and Nair, 2024; Xu et al., 2024; Xu et al., 2025; Zhang et al., 2025; Kim et al., 2025). More recently, reinforcement learning and deep reinforcement learning frameworks have shown strong scalability and robustness in planning and scheduling contexts (Shahab et al., 2024; Grumbach et al., 2024; Ngwu et al., 2025).

By contrast, quantum methods typically require additional and often cumbersome reformulations during the modelling stage, e.g. embedding linear constraints or multi-objective conditions into the QUBO representation through penalty terms, or introducing auxiliary variables to maintain feasibility (Rovara et al., 2024; Niu et al., 2025a). Such transformations inevitably lead to rapid growth in model size, which in turn restricts applicability and efficiency at realistic network scales (Cerezo et al., 2022). Moreover, recent quantum studies of VRP that tested QUBO encodings on actual or simulated quantum hardware explicitly indicate the performance bottlenecks in circuit depth, noise, and connectivity constraints, with most of the benchmark instances falling behind the strong classical baselines (Onah and Michielsen, 2025). These achievements suggest that, at present, quantum optimisation in transport remains largely exploratory rather than a demonstrably superior alternative to established classical methods.

Accordingly, this study adopts a rational and reflective perspective to explore the role of quantum optimisation in urban transport planning. While continued exploration of quantum methods may foster interdisciplinary exchange and contribute to longer-term methodological advances, it seems premature to portray them at this stage as a central solution for real-world urban transport problems. Therefore, this paper critically reviews the common optimisation problems in urban transport planning, highlights their effective solutions under classical algorithmic frameworks, as well as examines the modelling and application constraints associated with quantum methods. This study does not aim to dismiss quantum optimisation, but to clarify its current position as an exploratory complement rather than a foundational tool in urban transport planning.

**2. Quantum optimisation and classical methods**

*2.1 Overview of quantum methods*

Quadratic unconstrained binary optimisation (QUBO) has emerged as the most widely-adopted formulation framework in quantum optimisation. Its core principle is to absorb problem constraints into the objective function through penalty terms, thereby transforming the original model into a purely quadratic binary form. This formulation provides a unifying representation for a broad class of combinatorial problems. However, existing studies consistently report the challenges related to penalty weight calibration, model scalability, and feasibility preservation, especially when the problem size increases (Glover et al., 2022).

Quantum annealing (QA) exploits quantum tunnelling mechanisms to explore complex energy landscapes and has been tested on solving several transport-related optimisation problems, including traffic signal control, routing, and flow allocation. Most reported applications are small- or medium-scale instances. Restricted performance is continuously seen due to embedding overhead, hardware noise, and scalability constraints, and hybrid quantum–classical solvers deployed on current hardware platforms have not demonstrated consistent advantages over strong classical baselines in terms of either solution quality or computational efficiency (Quinton et al., 2025).

Variational quantum algorithms, most notably the Quantum Approximate Optimisation Algorithm (QAOA), have also attracted attention as candidates for near-term quantum devices. QAOA combines quantum circuits with classical optimisation loops and has shown promising prospects on selected structured instances. Nonetheless, recent reviews have highlighted the persistent challenges related to circuit depth, parameter trainability, and noise sensitivity, with empirical evidence of robust and repeatable advantages over classical methods remaining limited (Zhou et al., 2020; Abbas et al., 2024).

Overall, while these quantum optimisation methods introduce new computational perspectives, existing transport-related studies largely focus on benchmark testing and problem compatibility. A common structural constraint across these frameworks is their reliance on binary variable representations, which requires continuous decision variables to be discretised or augmented with auxiliary constructs. Such modelling transformations increase the problem dimensionality and further limit applicability in realistic urban transport planning contexts.

*2.2 Classical optimisation methods*

In the fields of transport and logistics, classical optimisation methods have developed into a rich and mature toolkit. These methods have shown dominant places in both academic research and practical applications in engineering and policy making. Depending on the level of precision and computational efficiency, they can be grouped into exact algorithms, heuristic and metaheuristic methods, and artificial intelligence-based optimisation methods.

Mixed integer linear programming (MILP) is the most widely-used modelling and solution framework. It can guarantee global optimality while handling complex factors such as capacity constraints, time windows and multi-objective formulations in transport networks. With the advances in commercial solvers such as CPLEX and Gurobi, MILP formulations of large-scale transport and logistics problems have become increasingly tractable (Liu et al., 2021; Li et al., 2022a; Tóth et al., 2024). At the same time, decomposition methods are critical for the improvement of computational efficiency. For example, Benders decomposition separates the original problem into master and subproblems in order to reduce the dimensionality of computation (Weng et al., 2024; Hosseini and Turner, 2024). Column generation shows effectiveness in dealing with exponentially large variable sets in vehicle routing and scheduling (Cóccola et al., 2015; Morabit et al., 2021). Branch and Bound provides a systematic exact solution procedure for integer optimisation problems (Yokoyama et al., 2015; Yokoyama et al., 2019; Adelgren and Gupte, 2022; Du et al., 2025). These methods have all demonstrated different extents of effectiveness in various transport applications.

The rise of artificial intelligence has introduced new methods to the optimisation problems in transport and logistics. Reinforcement learning (RL) and deep reinforcement learning (DRL) improve the optimisation performance through interactions with the environment and have been applied to traffic signal control, vehicle scheduling and route choice, with strong adaptability and robustness (Kumar et al., 2021; Haddad et al., 2022; Miletić et al., 2022; Yang and Fan, 2025; Karbasi and Yang, 2025). In addition, neural network-based optimisation methods can learn heuristic rules or approximate value functions, allowing them to generate near-optimal solutions in complex environments. These data-driven methods are well-suited to the uncertainty and dynamics in problems involving, thus extending the applicability of traditional optimisation frameworks.

Overall, classical optimisation methods in transport and logistics have developed into a multi-level toolkit. Specifically, exact algorithms ensure optimality, metaheuristics provide efficiency for large-scale problems, and artificial intelligence methods address dynamic and uncertain environments. These methods have already demonstrated a high level of maturity and practical value. By comparison, quantum optimisation should be regarded as an exploratory supplement rather than a core replacement at the current stage.

*2.3 Comparative perspective (classical vs. quantum methods)*

Overall, classical and quantum optimisation methods exhibit notable differences in their applicability to urban transport problems. Classical methods provide reliable and well-established solutions for medium- to large-scale networks. Exact algorithms supported by advanced solvers can address models with large numbers of variables, while heuristic and metaheuristic methods demonstrate efficient convergence in large-scale settings. In addition, artificial intelligence–based methods have shown strong adaptability and robustness in dynamic and uncertain transport environments.

By contrast, quantum optimisation methods remain constrained in practice by hardware scale, noise, and circuit depth despite their theoretical attraction due to mechanisms such as quantum tunnelling that may facilitate the exploration of complex solution landscapes. As a result, most existing studies are limited to small-scale or synthetic benchmark instances, which restricts the ability to demonstrate the consistent advantages on realistic urban transport networks. Moreover, quantum methods typically require additional reformulation steps, such as embedding linear constraints, multi-objective conditions, or logical relations into QUBO representations using penalty terms and binary variables. These transformations often lead to rapid growth in model size and increase the difficulty of maintaining feasibility. In comparison, classical optimisation frameworks can directly represent such constraints in a transparent and flexible manner.

Comparative studies to date generally indicate that classical methods remain more effective in terms of solution quality and computational efficiency in practical urban transport applications. At the same time, the primary contribution of quantum optimisation lies less in immediate performance gains and more in offering alternative modelling perspectives and problem-solving methods that may help broaden the optimisation research landscape. At the current stage of technological development, quantum optimisation is considered as a more appropriate exploratory complement to, rather than a replacement for, established classical methods.

## 3. Case comparison: typical optimisation problems in transport research

In the field of transport, optimisation is widespread with diverse forms, covering facility location, network optimisation, vehicle scheduling, network design and multi-objective disaster decision-making. To enable a systematic comparison, this paper selects five representative cases that together form a minimal yet sufficient set spanning strategic, tactical and operational levels, as shown in **Table 1**. The selection of these five cases is based on four criteria. First, structural representativeness. These cases cover static and dynamic settings, deterministic and uncertain conditions, single-objective and multi-objective formulations, as well as both resource allocation and routing decisions. Second, methodological comparability. Each of these cases has mature classical baselines such as MILP, decomposition, metaheuristics and reinforcement learning, along with documented attempts at quantum reformulation, which allows direct evaluation. Third, practical relevance. The selected cases frequently appear in policy making and engineering practice and have a significant decision-making impact. Fourth, diagnostic value. All the cases contain binary combinatorial structures and complex constraints that test the feasibility of QUBO transformations, scalability and constraint satisfaction in quantum optimisation. On this basis, this paper conducts a tabular comparison that summarises four aspects of problem category, classical methods, quantum reformulation, and critical evaluation. The aim is to help clearly identify the boundaries and differences between the two methods.

**Table 1.** Systematic comparison of classical and quantum methods for representative optimisation problems in transport research.

| Problem category | Classical methods (core modelling and algorithms) | Maturity and applications of classical methods | Quantum reformulation pathway | Practical limitations of quantum methods | Critical evaluation |
|---|---|---|---|---|---|
| Facility location and shelter allocation (Hansuwa et al., 2022; Eriskin and Karatas, 2024; Chacón-Tibaduiza et al., 2025) | MILP, Benders decomposition, heuristics and metaheuristics (GA, TS., et al.) | Able to solve medium and large-scale city and regional instances, with proven applications in emergency shelter and infrastructure planning | Reformulated as QUBO with capacity and coverage constraints embedded as penalty terms (Bhatia and Sood, 2024; Zhao et al., 2024) | Rapid growth of auxiliary variables, exponential expansion of QUBO size | Classical methods have already reached policy and engineering implementation; quantum remains demonstrative only |
| Route and network optimisation (Esposito Amideo et al., 2019; Wang et al., 2022) | Shortest path, traffic equilibrium models, operations research models, reinforcement learning (RL and DRL) | Real-time or near-real-time solutions for large-scale networks, capable of capturing congestion propagation and departure time choice | Routes and departure times binarised and reformulated as QUBO (Vinil et al., 2025; Dixit et al., 2024) | Large number of binary variables, difficult to extend to urban-scale networks | Classical methods are efficient and robust; quantum shows no substantive advantage |
| Vehicle routing and scheduling (VRP, CVRP, EVRP) (Eglese et al., 2006; Potvin et al., 2006; Maden et al., 2010) | Metaheuristics, MILP with decomposition | Widely applied in industry, capable of solving hundreds of vehicles and customers in practical instances | Small-scale VRP encoded into QUBO and tested on quantum annealers (Azad et al., 2023; Khalid et al., 2025) | Sensitive to noise, feasible only for toy problems with a few dozen nodes | Classical methods have reached industrial standards, and quantum has shown no scalable breakthrough |
| Network design and fairness optimisation (Caggiani et al., 2017; Zhang and Waller, 2019) | Benders decomposition, ε-constraint, multi-objective programming | Able to incorporate fairness objectives such as service balance and regional coverage, with mature applications in disaster resilience and infrastructure investment | Fairness constraints binarised with many auxiliary variables in QUBO (Dixit and Niu, 2023; Farahmand-Tabar and Afrasyabi, 2024) | Loss of structural properties and interpretability of fairness metrics | Classical methods balance efficiency and fairness; quantum modelling is unrealistic and distorts the essence |
| Disaster response and multi-objective optimisation (Ye et | Multi-objective evolutionary algorithms (such as NSGA- | Capable of handling uncertainty and generating | Simplified binary objectives reformulated as QUBO (Singh | Inability to capture uncertainty and difficulty in | Classical methods capture complexity and provide |
| al., 2024; Gupta et al., 2024; Bakhshian and Martinez-Pastor, 2024) | II), robust optimisation, scenario-based methods | Pareto frontiers, widely applied in emergency scheduling and disaster response | et al., 2025; Danach et al., 2025) | preserving multi-objective trade-offs | transparency; quantum oversimplification leads to distortion |

As shown in **Table 1**, classical optimisation methods have been systematically studied and widely applied in the transport field. These methods include exact algorithms, decomposition methods, and heuristic or intelligent techniques, all demonstrating mature and reliable performance across different scales and complexities. In contrast, quantum optimisation, despite its promising prospects, still faces many challenges, like modelling complexity, variable expansion, and hardware limitations when reformulated into QUBO, which restricts most outcomes to small-scale demonstrations and compatibility validation. Overall, quantum methods have yet to exhibit advantages over classical methods in real-world network scales, and their academic value lies more in exploratory complementarity rather than being an inevitable tool for practical applications.

## 4. Prospects and breakthrough pathways

Building on the preceding analysis and comparison of the structural limitations of quantum optimisation in transport systems, this section further discusses the appropriate role and potential development paths of quantum optimisation in the transport field under current technological conditions. It should be emphasised that there is no specific technical roadmap proposed in this paper. Instead, from the perspective of research methods and problem matching, we examine how to avoid technology-driven misjudgements and provide more rational guidance for future research.

*4.1 Redefining the role of quantum optimisation in urban transport research*

Existing evidence shows that quantum hardware and related methods do not yet have the capacity to act as dominant tools that can fully replace classical methods. In transport research, quantum optimisation is normally used as a complementary tool to support specific types of subproblems or analytical stages, rather than to finish full system-wide tasks. This is particularly relevant for the problems involving large-scale continuous variables, multiple constraints, and complex decision structures, where classical methods still offer obvious advantages in modelling transparency, system coherence, and interpretability.

Based on this, **Table 2** provides a conceptual overview of the appropriate positioning of classical optimisation methods and quantum optimisation across typical transport problems, considering different decision tiers and problem characteristics. This table clarifies the scope of applicability and the potential division of roles between different optimisation methods in the field of transport research.

**Table 2.** Conceptual positioning of quantum optimisation in transport research.

| Problem characteristics | Typical examples | Core strengths of classical methods | Potential role of quantum optimisation | Key current limitations |
|---|---|---|---|---|
| Strategic-level planning | Network design, long-term investment, and infrastructure layout | Global consistency, multi-objective trade-offs, strong interpretability | Not suitable as a primary solution approach | Binary reformulation distorts the structure, and scalability is unmanageable |
| Tactical-level decisions | Resource allocation, scheme screening, scenario comparison | Stable solution quality, mature modelling, policy relevance | Local combinatorial search, exploratory scenario analysis | QUBO expansion, complex constraint embedding |
| Operational-level problems | Scheduling, routing, short-term control | Real-time capability, robustness, and operational stability | Small-scale subproblem exploration | Hardware limitations and latency constraints |
| Problems dominated by continuous variables | Traffic assignment, time-dependent demand modelling | Direct representation without discretisation | Forced discretisation required | Precision loss and rapid variable growth |
| Highly discrete but scale-limited problems | Combinatorial selection, substructure matching | Constrained search space | Potential exploratory value | Noise sensitivity and solution stability |
| Multi-objective and uncertainty-driven | Efficiency-equity balance, disaster | Pareto frontier generation | Difficult to preserve trade-off structure | Objective simplification and information loss |

| problems | response planning | | |
|---|---|---|---|

*4.2 Shifting from technology-driven to problem-driven research pathways*

This study argues that future research should not focus on continuously searching for transport problems that fit quantum computing frameworks. Instead, attention shall be paid to which transport problems may genuinely benefit from quantum technologies at a structural level. Rather than pursuing global optimal solutions at the full system scale, quantum optimisation at the current stage is more likely to show exploratory value in the scenarios that involve highly-discrete decisions, high combinatorial complexity, and limited problem scale. Such problems usually emphasise multi-option screening, extreme scenario identification, or local search acceleration, rather than one-time solution of complete planning tasks. On this basis, **Table 3** explicitly identifies the transport problems that quantum technologies may reasonably handle. This table summarises the application scenarios that show potential exploratory value from a problem structure perspective, and conceptually outlines the appropriate ways to use quantum optimisation in terms of structural features, potential mechanisms, and research role positioning.

**Table 3.** Problem-driven opportunities for exploratory quantum optimisation in transport planning.

| Specific transport planning problem | Structural characteristics | Why may quantum optimisation have exploratory value? | Appropriate research role |
|---|---|---|---|
| Identification of critical links and nodes (network vulnerability analysis) | Binary selection, combinatorial growth with candidate set size | Combinatorial effects of multiple critical elements are difficult to enumerate | Exploratory identification of high-risk link or node combinations |
| Traffic congestion scenario analysis | Discrete congestion state combinations, focus on extreme or high-impact cases | Large scenario space with emphasis on worst-case configurations | Exploratory analysis and stress testing of congestion scenarios |
| Failure combination screening under disaster scenarios | Discrete failure combinations across links or facilities | Exhaustive enumeration is infeasible; focus on extreme disruption patterns | Screening of critical failure scenarios |
| Path planning and route selection | Discrete path choices within constrained subnetworks | Combinatorial growth of feasible paths under constraints | Exploratory search over limited path sets or subnetwork routing |
| Traffic signal control scheme combination optimisation | Discrete phase and timing combinations across intersections | High combinatorial complexity in multi-intersection coordination | Exploratory optimisation of small-scale subproblems |
| Facility location problems (e.g. shelters, charging stations) | Binary siting decisions with explicit coverage or capacity constraints | Limited but combinatorial option sets are suitable for pre-screening | Generation of candidate location schemes |
| Transport infrastructure investment option selection | Discrete investment choices with limited alternatives | Comparison of multiple option combinations rather than exact optima | Preliminary screening of investment options |

**4.3 Hybrid frameworks as a feasible transitional pathway**

In the foreseeable future, quantum optimisation has great potential to be embedded into existing transport systems through hybrid frameworks. Under such frameworks, classical methods handle continuous decision modelling, complex constraint management, and system coherence, while quantum methods focus on local discrete combinatorial search or multi-solution exploration tasks. This division of roles arises from the differences in their problem representation capability and computational. Based on current hardware and algorithms, hybrid frameworks already show a relatively clear technical pathway. For example, quantum annealing hardware platforms such as the D-Wave series usually represent combinatorial optimisation problems using binary or discrete variable models, including binary quadratic model (BQM), constrained quadratic model (CQM), and discrete quadratic model (DQM). These models provide flexibility for small-scale binary choices, limited discrete decisions, or combinatorial screening subproblems (Niu et al., 2025a). However, when facing large-scale continuous variables or complex constraint structures, classical methods are still the preferred choice for external coordination and constraint management. In parallel, gate-based variational quantum algorithms such as QAOA are often viewed as potential methods for near-term quantum devices. Their application contexts also concentrate on scale-limited and structurally well-defined discrete combinatorial problems. In transport research, such algorithms are more suitable for exploratory analysis or subproblem solving than classical methods. **Table 4** specifies the problem types that suit different quantum optimisation models and algorithms.

**Table 4.** Conceptual mapping between quantum optimisation models and their roles in hybrid transport optimisation frameworks.

| Quantum model or algorithm | Typical hardware or approaches | Transport subproblems suitable for embedding | Contexts where the model is not appropriate |
| --- | --- | --- | --- |
| BQM | Quantum annealing (e.g. D-Wave) | Binary combinatorial selection, screening of critical elements | Problems dominated by continuous variables, large-scale network optimisation |
| CQM | Quantum annealing (e.g. D-Wave) | Discrete decision subproblems with a limited number of linear constraints | System-level planning problems with complex and dense constraint structures |
| DQM | Quantum annealing (e.g. D-Wave) | Selection among multiple discrete options with a bounded problem size | High-dimensional continuous decisions or real-time traffic control |
| QAOA | Gate-based quantum computing | Structured and size-limited discrete combinatorial problems | Full-network route planning and dynamic traffic management |

*4.4 Implications for transport research and policy practice*

Based on existing findings, the appropriate role of quantum optimisation in transport research and policy making can be summarised in three points.

(1) Future studies should not treat the achievement of quantum advantages as the main evaluation criterion. Research should instead focus on how well quantum methods fit different transport decision stages. Clarification should be conducted about which problems may benefit from new computational methods due to their structure, and which problems remain better handled by the classical methods.

(2) Before quantum technology shows stable and reproducible practical benefits, public investment should remain cautious and should not position it as a core tool for solving real-world transport problems. A more suitable strategy is to support method validation, scenario matching, and prototype exploration.

(3) The potential value of quantum computing in transport mainly lies in how it expands problem modelling methods and analytical perspectives. For this reason, coordinated exploration across transport planning, operational optimisation, computational science, and policy research carries greater long-term value than a narrow focus on technical breakthroughs alone.

## 5. Conclusion

This study conducts a systematic investigation on the application of quantum optimisation in transport systems, mainly focused on suitability, limitations, and practicality. Rather than following technology-driven narratives centred on computational performance or potential quantum advantage, this study stresses the need to examine the appropriate role of emerging technologies from the perspective of problem structure and decision tiers.

It has been obviously shown that classical optimisation methods have reached a high level of maturity at strategic, tactical, and operational tiers. They demonstrate clear advantages in modelling transparency, scalability, and engineering feasibility, and have been widely applied in real-world transport systems. By contrast, quantum optimisation methods still face multiple constraints during the modelling process, including variable scale expansion, limited constraint representation, and insufficient hardware capability. Their application scope remains largely confined to small-scale problem validation or exploratory studies.

On this basis, this study further argues that, in the foreseeable future, quantum optimisation does not necessarily serve as a core tool to replace existing transport planning methods. It can act as a complementary approach within hybrid frameworks. In this way, quantum optimisation can support specific tasks such as combinatorial search, scenario screening, or solution exploration. That means, classical methods maintain system coherence and handle continuous decision modelling, while quantum methods address discrete subproblems with constrained structure and controllable scale.

Overall, by clarifying both the application scenarios and practical limits of quantum optimisation, this study aims to provide exact and sustainable guidance for future transport research and policy making. It seeks to promote responsible exploration and well-grounded application of emerging computational technologies in urban transport systems.

**Declaration of generative AI and AI-assisted technologies in the writing process**
During the preparation of this work, we employed ChatGPT-5.2 solely to assist with language editing and polishing. No content was generated by AI. After using this tool, we thoroughly reviewed and edited the content as necessary and accept full responsibility for the content of the published article.

**Appendix A:** Representative literature review on quantum computing in transport research

| Authors | Data | Title | Journal |
| --- | --- | --- | --- |
| Neukart F, Compostella G, Seidel C, Von Dollen D, Yarkoni S, Parney B | 2017 | Traffic flow optimization using a quantum annealer | *Frontiers in ICT* |
| Feld, S., Roch, C., Gabor, T., Seidel, C., Neukart, F., Galter, I., Mauerer, W. and Linnhoff-Popien, C | 2019 | A hybrid solution method for the capacitated vehicle routing problem using a quantum annealer | *Frontiers in ICT* |
| Stollenwerk, T., O'Gorman, B., Venturelli, D., Mandra, S., Rodionova, O., Ng, H., Sridhar, B., Rieffel, E.G. and Biswas, R | 2019 | Quantum annealing applied to de-conflicting optimal trajectories for air traffic management | *IEEE transactions on intelligent transportation systems* |
| Hussain, H., Javaid, M.B., Khan, F.S., Dalal, A. and Khalique, A | 2019 | Optimal control of traffic signals using quantum annealing | *arXiv preprint* |
| Inoue, D., Okada, A., Matsumori, T., Aihara, K. and Yoshida, H | 2021 | Traffic signal optimization on a square lattice with quantum annealing | *Scientific reports* |
| Dixit, V.V. and Niu, C | 2023 | Quantum computing for transport network design problems | *Scientific reports* |
| Xu, H.Z., Chen, J.H., Zhang, X.C., Lu, T.E., Gao, T.Z., Wen, K. and Ma, Y | 2023 | High-speed train timetable optimization based on space–time network model and quantum simulator | *Quantum Information Processing* |
| Dixit, V.V., Niu, C., Rey, D., Waller, S.T. and Levin, M.W | 2024 | Quantum computing to solve scenario-based stochastic time-dependent shortest path routing | *Transportation Letters* |
| Niu, C., Rastogi, P., Soman, J., Tamuli, K. and Dixit, V | 2025 | Applying Quantum Computing to Solve Multicommodity Network Flow Problem | *IEEE Transactions on Intelligent Transportation Systems* |
| Salloum, H., Zhanalin, S., Badr, A.A. and Kholodov, Y | 2025 | Mini-scale traffic flow optimization: an iterative QUBOs approach converting from hybrid solver to pure quantum processing unit | *Scientific reports* |
| Koniorczyk, M., Krawiec, K., Botelho, L., Bešinović, N. and Domino, K | 2025 | Solving rescheduling problems in heterogeneous urban railway networks using hybrid quantum-classical approach | *Journal of Rail Transport Planning & Management* |
| Du, Z., Wandelt, S. and Sun, X | 2026 | Overcoming computational challenges in air transportation: A quantum computing perspective of the status quo and future applicability | *Transportation Research Part C: Emerging Technologies* |